\documentclass[leqno, fleqn]{amsart}%
\usepackage{amsmath}
\usepackage{amsfonts}
\usepackage{amssymb}%

\makeatletter \DeclareMathSymbol{\Gamma}{\mathalpha}{letters}{"00}
\DeclareMathSymbol{\Delta}{\mathalpha}{letters}{"01}
\DeclareMathSymbol{\Theta}{\mathalpha}{letters}{"02}
\DeclareMathSymbol{\Lambda}{\mathalpha}{letters}{"03}
\DeclareMathSymbol{\Xi}{\mathalpha}{letters}{"04}
\DeclareMathSymbol{\Pi}{\mathalpha}{letters}{"05}
\DeclareMathSymbol{\Sigma}{\mathalpha}{letters}{"06}
\DeclareMathSymbol{\Upsilon}{\mathalpha}{letters}{"07}
\DeclareMathSymbol{\Phi}{\mathalpha}{letters}{"08}
\DeclareMathSymbol{\Psi}{\mathalpha}{letters}{"09}
\DeclareMathSymbol{\Omega}{\mathalpha}{letters}{"0A}
\DeclareMathSymbol{\varGamma}{\mathalpha}{operators}{"00}
\DeclareMathSymbol{\varDelta}{\mathalpha}{operators}{"01}
\DeclareMathSymbol{\varTheta}{\mathalpha}{operators}{"02}
\DeclareMathSymbol{\varLambda}{\mathalpha}{operators}{"03}
\DeclareMathSymbol{\varXi}{\mathalpha}{operators}{"04}
\DeclareMathSymbol{\varPi}{\mathalpha}{operators}{"05}
\DeclareMathSymbol{\varSigma}{\mathalpha}{operators}{"06}
\DeclareMathSymbol{\varUpsilon}{\mathalpha}{operators}{"07}
\DeclareMathSymbol{\varPhi}{\mathalpha}{operators}{"08}
\DeclareMathSymbol{\varPsi}{\mathalpha}{operators}{"09}
\DeclareMathSymbol{\varOmega}{\mathalpha}{operators}{"0A}

\newcommand{\allmodesymb}[2]{\relax\ifmmode{\mathchoice
{\mbox{\fontsize{\tf@size}{\tf@size}#1{#2}}}
{\mbox{\fontsize{\tf@size}{\tf@size}#1{#2}}}
{\mbox{\fontsize{\sf@size}{\sf@size}#1{#2}}}
{\mbox{\fontsize{\ssf@size}{\ssf@size}#1{#2}}}} \else
\mbox{#1{#2}}\fi}

\makeatother
\makeatletter
\renewcommand*\subjclass[2][2000]{%
  \def\@subjclass{#2}%
  \@ifundefined{subjclassname@#1}{%
    \ClassWarning{\@classname}{Unknown edition (#1) of Mathematics%
      Subject Classification; using '2000'.}%
  }{%
    \@xp\let\@xp\subjclassname\csname subjclassname@#1\endcsname%
  }%
} \makeatother

\theoremstyle{plain}

\theoremstyle{remark}

\allowdisplaybreaks
\numberwithin{equation}{section}
\begin{document}
\title[Complete K\"{a}hler--Einstein metric on $Y_{II}(r,p;K)$]
{Complete Einstein-K\"{a}hler Metric and Holomorphic Sectional
Curvature on $Y_{II}(r,p;K)$}
\date{10th June 2005}

\author{Weiping YIN }
\address{W YIN: Dept. of Math., Capital Normal Univ., Beijing 100037, China}
\email{wyin@mail.cnu.edu.cn; wpyin@263.net}
\author{Liyou ZHANG }
\address{L ZHANG: Dept. of Math., Capital Normal Univ., Beijing
100037, China} \email{zhangly@mail.cnu.edu.cn} \subjclass{32H15,
32F07, 32F15} 
\thanks{Project supported in part by NSF of China (Grant NO.
10471097) and the Doctoral Programme Foundation of NEM of China.}

\begin{abstract}
The explicit complete Einstein-K\"{a}hler metric on the second
type Cartan-Hartogs domain $Y_{II}(r,p;K)$ is obtained in this
paper when the parameter $K$ equals $\frac p2+\frac 1{p+1}$. The
estimate of holomorphic sectional curvature under this metric is
also given which intervenes between $-2K$ and $-\frac{2K}p$ and it
is a sharp estimate. In the meantime we also prove that the
complete Einstein-K\"ahler metric is equivalent to the Bergman
metric on $Y_{II}(r,p;K)$ when $K=\frac p2+\frac 1{p+1}$.

\end{abstract}

\maketitle

\section*{Introduction}

It is well known that the Bergman, Carath$\acute{\rm e}$odory,
Kobayashi and Einstein-K$\ddot{\rm a}$hler metrics are four
classical invariant metrics in complex analysis. The Bergman
metric was introduced by S.Bergman for one variable $^{[1]}$in
1921 and for several variables$^{[2]}$ in 1933.C.Carath$\acute{\rm
e}$odory introduced the invariant distance in 1926$^{[3]}$ and
H.Reiffen introduced the invariant metric in 1963$^{[4]}$,
therefore the Carath$\acute{\rm e}$odory metric is also called
Carath$\acute{\rm e}$odory-Reiffen metric. The Kobayashi metric
was introduced by S.Kobayashi in 1967$^{[5]}$ and by H.Royden in
1970$^{[6]}$. Therefore the Kobayashi metric is also called
Kobayashi-Royden metric. Let $M$ be a complex manifold, then a
Hermitian metric $\sum_{i,j}g_{i,\overline{j}}$d$z^i\bigotimes
$d$\overline{z}^j$ defined on $M$ is said to be K$\ddot{\rm
a}$hler if the K$\ddot{\rm a}$hler form
$\sl\Omega=\sqrt{-1}\sum_{i,j}g_{i,\overline{j}}$d$z^i\bigwedge$
d$\overline{z}^j$ is closed. The Ricci form of this metric is
defined to be $-\partial \overline{\partial}{\rm{log}}
{\rm{det}}(g_{i,\overline j})$. If the Ricci form of the K\"ahler
metric is proportional to the K\"ahler form, the metric is called
Einstein-K\"ahler. If the manifold is not compact, it requires the
metric to be complete.According to a famous article$^{[7]}$ of Wu,
one knows that Einstein-K\"ahler metric is the most difficult to
compute among the four metrics because its existence is proved by
complicate nonconstructive methods. If we normalize the metric by
requiring the scalar curvature to be minus one, then the
Einstein-K\"ahler metric is unique.

Cheng and Yau$^{[8]}$ proved that any bounded pseudo convex domain
$D$ with continuous second partial derivatives boundary admits a
complete Einstein-K\"ahler metric. Mok and Yau have extended this
result to an arbitrary bounded pseudo convex domain in
$\mathbb{C}^{n}$$^{[9]}$. If this Einstein-K\"ahler metric is
given by
$$E_{D}(z):=\sum \frac{\partial^2g}{\partial z_i \partial
\overline{z}_j} dz_i \overline {dz_j}\ ,$$  then $g$ is the unique
solution to the boundary problem of the Monge-Amp\` ere equation:
$$
 \left\{
\begin{array}{ll}
\displaystyle {\rm{det}} \left( \frac{\partial^2g}{\partial z_i
\partial \overline{z}_j} \right) =e^{(n+1)g} & z \in  D, \\
\displaystyle g=\infty & z \in \partial D,\\
\end{array}
\right. $$ and $g$ is called generating function of $E_{D}(z).$
Obviously, if one obtains $g$ in explicit formula, then the
Einstein-K\"ahler metric is also explicit. We knew by far that the
explicit formulas for the Einstein-K\"ahler metric, however, are
only known on homogeneous domains, and it is exactly the Bergman
metric. We knew little as far as the nonhomogeneous domains
concerned.

In this paper, we just consider a class of nonhomogeneous domains
introduced by prof. W.Yin$^{[10]}$ and G.Roos in 1998. It is
defined as:
$$
Y_{II}(r,p;K)=\{w \in {\mathbb {C}^r},Z \in
R_{II}(p):|w|^{2}<\det(I-Z\overline Z)^{\frac{1}{K}} ,K>0
\}:=Y_{II},$$ where $R_{II}(p)$ is the second type of symmetric
classical domain and $\det$ is the usual determinant of matrices.
$p>1$ is a positive integer. $w$ is a vector with $r$ entries and
$|w|^2=|w_1|^2+|w_2|^2+\cdots+|w_r|^2$. Its Bergman kernel
function is given in explicit formula$^{[11]}$ and hence $Y_{II}$
is Bergman exhaustion, therefore, $Y_{II}$ is a bounded pseudo
convex domain which admits a unique complete Einstein-K\"ahler
metric.

We know that the explicit complete Einstein-K\"ahler metric on
$Y_{II}(1,p;K)$ has been obtained$^{[12]}$ when the parameter
satisfies $K=\frac p2+\frac 1{p+1}$. The corresponding generating
function $g$ has explicit form as follows:
$$\begin{array}{lll}
g&=&\frac{1}{n+1}\log[{Y}^{n+1}\det(I-Z\overline Z)^{-(1+p+\frac{1}{K})}K^{1-n}]\\
&=&\log[\frac{1}{1-X}\det(I-Z\overline
Z)^{-\frac{1}{K}}K^{\frac{1-n}{1+n}}]\ ,
\end{array}
$$
where the parameters $X$ and $Y$ are $X=|w|^2[\det(I-Z\overline
Z)]^{-\frac 1K}, Y=(1-X)^{-1}$ and $n=\frac{p(p+1)}{2}+1$ is the
dimension of $Y_{II}(1,p;K).$

From this result, we guess when $r>1,$ $K=\frac p2+\frac 1{p+1}$,
the generating function $g$ of complete Einstein-K\"ahler metric
on nonhomogeneous domain $Y_{II}(r,p;K)$ has the following form:
$$\begin{array}{lll}
g&=&\frac{1}{N+1}\log[{Y}^{N+1}\det(I-Z\overline Z)^{-(1+p+\frac{r}{K})}K^{r-N}]\\
&=&\log[\frac{1}{1-X} \det(I-Z\overline Z)^{-\frac
1K}K^{\frac{r-N}{1+N}}]\ ,
\end{array}
$$
where
\begin{align*}
X  &  =|w|^{2}[\det(I-Z\overline{Z})]^{-\frac{1}{K}}=(|w_{1}|^{2}+|w_{2}%
|^{2}+\cdots+|w_{r}|^{2})\det(I-Z\overline{Z})^{-\frac{1}{K}},\\
Y  &  =(1-X)^{-1}%
\end{align*}
and $N=\frac{p(p+1)}{2}+r$ is the dimension of $Y_{II}(r,p;K).$ We
will prove our conjecture correct through verifying $g$ is the
unique solution to the boundary problem of the Monge-Amp\` ere
equation.

This paper is organized as follows. Firstly, we present some basal
facts and results which we need about $Y_{II}(r,p;K)$. Secondly,
we will verify $g$ satisfies Monge-Amp\` ere equation and the
Direchlet boundary condition, so it gives the complete
Einstein-K\"ahler metric. Thirdly, an estimate of the holomorphic
sectional curvature under this metric is given. Finally, we prove
that the complete Einstein-K\"ahler metric is equivalent to the
Bergman metric on $Y_{II}(r,p;K)$ in the case $K=\frac p2+\frac
1{p+1}$.

\section{Preliminaries }

In this section, we give a few lemmas about $Y_{II}(r,p;K)$ which
will be needed later.\\
\\
 {\bf Lemma 1.} $\mathrm {Aut}(Y_{II})$ indicates
the holomorphic automorphism group of $Y_{II}(r,p;K)$ consisting
of the following mappings:
$$ \left\{
     \begin{array}{lll}
         w_j^*& =& w_j\det{(I-Z_0\overline Z_0)^{\frac1{2K}}}
         \det{(I-Z\overline Z_0)^{-\frac1K}}, \ \ \ j=1,2,\cdots,r. \\
         Z^*& =& A(Z-Z_0){(I-{\overline
         Z_0}Z)^{-1}}{{\overline{A}}\!^{-1}}.
         \end{array}
      \right.$$
where $\overline{A}^t$ denotes the conjugation and transpose of
$A$ and $\overline{A}^tA=(I-Z_0\overline Z_0)^{-1}, {Z_0} \in
{R_{II}(p)}.$
\\
{\bf Proof.} See ref. $[11]$.
\\

Obviously, any of the above mappings maps the point $(w,Z_0)$ onto
the point $(w^*,0)$ and
$Z^*=A(Z-Z_0)(I-\overline{Z_0}Z)^{-1}\overline{A}^{-1}$ is the
holomorphic automorphism of $R_{II}(p)$.\\
\\
 {\bf Lemma 2.} Let $X=X(Z,w)=|w|^2[\det(I-Z\overline
Z)]^{-1/K}$, then $X$ is invariant under the mapping of
$\mathrm{Aut}(Y_{II})$. That is
$X(Z^*,w^*)=X(Z,w)$.\\
 {\bf Proof.} See ref. $[11]$.
\\

Hereafter  we write $Z\in R_{II}(p)$ as $$
   Z=Z^t=\left(
    \begin{array}{cccc}
   z_{11}&\displaystyle \frac 1{\sqrt2}z_{12}&\cdots& \displaystyle \frac 1{\sqrt2}z_{1p}\\
   \displaystyle \frac 1{\sqrt2}z_{21}&z_{22}&\cdots& \displaystyle \frac 1{\sqrt2}z_{2p}\\
   \displaystyle \vdots&\vdots&\ddots&\vdots\\
   \displaystyle \frac 1{\sqrt2}z_{p1}& \displaystyle \frac 1{\sqrt2}z_{p2}&\cdots&z_{pp}
    \end{array}
    \right).$$
    and$$
z=(z_{11},z_{12},\cdots,z_{1p},z_{22},z_{23},\cdots,z_{2p},\cdots
\cdots,z_{pp})$$ is the $1\times p(p+1)/2$ matrix. ${Z}^t$ denotes
the transpose of $Z$.\\
\\
{\bf Lemma 3.} Suppose $(Z^*,w^*)=F(Z,w)\in \mathrm{Aut}(Y_{II})$
which maps $(Z_0,w)$ onto $(0,w^*)$; let J$_F$ be the Jacobi
matrix of $F(Z,w)$, i.e.
$$ \mathrm{J}_F=\left(
\begin{array}{cc}
\displaystyle \frac{\partial z^*}{\partial z}& \displaystyle
\frac{\partial w^*}{\partial z}\\
\\
\displaystyle 0& \displaystyle \frac{\partial w^*}{\partial w}
\end{array}
\right). $$ Then one has
$$\begin{array}{l} \displaystyle{\bigg(\frac{\partial z^*}{\partial z}\bigg)}_{Z_0=Z} =[{A^\prime}\cdot\!\!\times {A^\prime}]_s,\\
\\
\displaystyle {\bigg( \frac{\partial w^*}{\partial z}\bigg)}_{Z_0=Z}=\frac{1}{K} \det(I-Z\overline Z)^{-\frac{1}{2K}}E(Z)^tw,\\
\\
\displaystyle{\bigg( \frac{\partial w^*}{\partial
w}\bigg)}_{Z_0=Z}=I\det(I-Z\overline Z)^{-\frac{1}{2K}},
\\
\end{array}
$$
where $E(Z)$ is the $1\times{p(p+1)/2}$ matrix
$$E(Z)=(tr[(I-Z\overline Z)^{-1} I^*_{11}\overline Z],
\cdots,tr[(I-Z\overline Z)^{-1} I^*_{kl}\overline Z],\cdots,
tr[(I-Z\overline Z)^{-1} I^*_{pp}\overline Z])$$  and
$$ I^*_{kl}= \left\{
  \begin {array}{ll}
  \displaystyle {\frac 1{\sqrt2}}(I_{kl}+I_{lk}), &\quad k<l, \\
  I_{kk},&\quad k=l.
  \end{array}
  \right.
$$
Here $I_{kl}$ is defined as a $p\times p$ matrix, the $(k,l)$th
entry of $I_{kl}$, i.e. the entry located at the junction of the
$k$-th row and $l$-th column of $I_{kl}$ is 1, and others entries
of $I_{kl}$ are zero. The meaning of $ [A \cdot\!\!\times A]_s $
can be found in [13].
\\
{\bf Proof.} It can be got by direct computation.
\\
\\
{\bf Lemma 4.} If $(Z^*,w^*)=F(Z,w)\in \mathrm{Aut}(Y_{II})$ and
$T=T[(Z,w),\overline{(Z,w)}]$ is the metric matrix of  the
Einstein-K\"ahler metric on $Y_{II}(r,p;K)$, then one has
$$ T[(Z,w),\overline{(Z,w)}] =[\mathrm{J}_FT[(Z^*,w^*),\overline{(Z^*,w^*)}]
\overline{\mathrm{J}}_F^t ]_{Z_0=Z}, $$ and
$|\mathrm{J}_F|^2_{Z_0=Z}=\det(I-Z\overline{Z})^{-(p+1+\frac{r}{K})}$,
where $|\mathrm{J}_F|=\det \mathrm{J}_F.$
\\
{\bf Proof.} It can be proved by using the invariance of the
Einstein-K$\ddot{\rm a}$hler metric under the holomorphic
automorphism of $Y_{II}$.
\\
\\
{\bf Lemma 5.} Let $Z$ be a $p\times p$ symmetric matrix, then the
following inequality holds:
$$ tr(Z\overline{Z}Z\overline{Z})\leq tr(Z\overline{Z})
tr(Z\overline{Z}) \leq p [tr(Z\overline{Z} Z\overline{Z})].$$
\\
 {\bf Proof.} The lemma is trivial when $p=1$, so we
consider the case $p>1$ only. Let $Z$ be a non-zero matrix, there
exists a unitary matrix $U$ such that
$$ Z=U^t\left(
\begin{array}{ccc}
\lambda_1&0 &0 \\ 0&\ddots&0\\ 0&0&\lambda_p\\
\end{array}
\right )  U ,\quad \lambda_1\geq \lambda_2\geq \cdots
\lambda_p\geq 0,\quad \lambda_1>0. $$ then  $$tr(Z\overline{Z}
Z\overline{Z})= \lambda_1^4+\lambda_2^4+\cdots +\lambda_p ^4\
,\quad tr(Z\overline{Z})= \lambda_1^2+\lambda_2^2+\cdots
+\lambda_p^2\ .$$ By using Cauchy-Schwartz inequality we have
$$
\begin{array}{lll}
\lambda_{1}^{4}+\cdots+\lambda_{p}^{4}  &
\leq(\lambda_{1}^{2}+\cdots
+\lambda_{p}^{2})^{2}=|(1^{2},\cdots,1^{2})\cdot(\lambda_{1}^{2}%
,\cdots,\lambda_{p}^{2})^{t}|^{2}\\
&  \leq{|(1^{2},\cdots,1^{2})|^{2}}{|(\lambda_{1}^{2},\cdots,\lambda_{p}%
^{2})|^{2}=}p(\lambda_{1}^{4}+\cdots+\lambda_{p}^{4}),
\end{array}
$$
i.e. $ tr(Z\overline{Z}Z\overline{Z})\leq tr(Z\overline{Z})
tr(Z\overline{Z}) \leq p [tr(Z\overline{Z} Z\overline{Z})].$

It is obvious that $tr(Z\overline{Z}Z\overline{Z})=
tr(Z\overline{Z}) tr(Z\overline{Z})$ holds if
$\lambda_2=\cdots=\lambda_p=0$ and $tr(Z\overline{Z})
tr(Z\overline{Z})= p [tr(Z\overline{Z} Z\overline{Z})]$ holds if
and only if $\lambda_1=\lambda_2=\cdots=\lambda_p$.

\section{Complete Einstein-K\"ahler metric with explicit formula}

 Let $Z\in R_{II}(p)$ and $(Z,w)\in Y_{II}(r,p;K)$. Denote
 $(z,w)=(z_{11},z_{12},\cdots,z_{1p},z_{22},$
 $\cdots,z_{2p},\cdots,z_{pp},w_1,\cdots,w_r)= (z_1,z_2, \cdots
,z_N),$ where $N=\frac{p(p+1)}2+r$ is the dimension of
$Y_{II}(r,p;K)$.

Note that $$g_{\alpha \overline \beta}(z,w)= \frac{\partial^2g}
{\partial z_\alpha \partial \overline{z}_{\beta}}, \qquad \alpha ,
\beta =1,2,\cdots,N,$$ where
$$\frac{\partial g}{\partial z_{p(p+1)/2+j}}=
\frac{\partial g}{\partial w_j}, \qquad j=1,\cdots,r. $$

In order to prove our conjecture is right, we have to verify that
the generating function $g$ is the unique solution to the
Dirichlet boundary problem of complex Monge-Amp\` ere equation:
$$ \left
\{
\begin{array}{ll}
  \det (g_{\alpha \overline \beta}(z,w))=e^{(N+1)g(z,w)}, &(z,w) \in Y_{II}, \\
 g=\infty, & (z,w) \in \partial Y_{II}.\\
\end{array}
\right. \eqno{(1)} $$

Let $F\in \mathrm{Aut}(Y_{II}), F(Z,w)=F(Z^*,w^*).$ According to
lemma 4 we know $$\det(g_{\alpha \overline \beta}(z,w))=
|\mathrm{J}_F|^2\det(g_{\alpha \overline \beta}(z^*,w^*)), $$
especially, if choose $Z_0=Z$ and denote holomorphic automorphism
${F|_{Z_0=Z}}$ by $F_0$, then
$$\det(g_{\alpha\overline \beta}(z,w))=  |\mathrm{J}_{F_0}|^2\det(g_{\alpha
\overline\beta}(0,w^*)).\eqno{(2)} $$ Hence we need only to know
the value of $\det(g_{\alpha \overline \beta}(z^*,w^*))$ at the
point $(0,w^*)$. We substitute $w$ for $w^*$ in the following
computation.

Since $Y=(1-X)^{-1}, X=|w|^2[\det(I-Z\overline Z)]^{-1/K},$ the
generating function $g$ can be rewritten by
$$ g=\log Y + \log X -\log |w|^2 + \frac{r-N}{1+N} \log K .$$
Thus
$$\begin{array}{lll}
\displaystyle\frac{\partial g}{\partial z_\alpha}&=&
\displaystyle\left(Y+X^{-1}\right)\frac{\partial X}{\partial
z_\alpha}\ ,\qquad \alpha=1,2,\cdots,N-r.\\
\displaystyle\frac{\partial g}{\partial w_{i}}&=&\displaystyle
Y\frac{\partial X}{\partial w_{i}}\ , \qquad i=1,2,\cdots,r.
\end{array}$$
and
$$
\begin{array}{lll}
\displaystyle \frac{\partial^2g}{\partial z_\alpha \partial
\overline{z}_\beta}&=& \displaystyle
\left(Y+X^{-1}\right)\frac{\partial^2 X}{\partial z_\alpha
\partial \overline{z}_\beta}
 + \left(Y^2-X^{-2}\right)
\frac{\partial X}{\partial z_\alpha}\frac{\partial X} {\partial \overline{z}_\beta}\ ,\\
\displaystyle \frac{\partial^2g}{\partial z_\alpha \partial
\overline{w}_j}&= &\displaystyle (Y+X^{-1})\frac{\partial^2
X}{\partial z_\alpha \partial \overline{w}_j}
 + (Y^2-X^{-2}) \frac{\partial X}{\partial z_\alpha}\frac{\partial X} {\partial \overline{w}_j}\ ,
\qquad j=1,2,\cdots,r\ ;\\
\displaystyle \frac{\partial^2 g}{\partial
w_{i}\partial\overline{z}_\beta} &= &\displaystyle Y
\frac{\partial^2 X}{\partial w_i\partial\overline{z}_\beta} + Y^2
\frac{\partial X}{\partial w_{i}}\frac{\partial X}{\partial
\overline{z}_\beta}\ ,
\qquad i=1,2,\cdots,r\ ;\\
\displaystyle \frac{\partial^2g}{\partial w_i
\partial\overline{w}_j}&=& \displaystyle Y \frac{\partial^2
X}{\partial w_i\partial \overline{w}_{j}}
 + Y^2\frac{\partial X}{\partial
w_{i}}\frac{\partial X}{\partial \overline{w}_{j}}\ ,\qquad
i,j=1,2,\cdots,r\ .
\end{array}\eqno{(3)}
$$
According ref.[12] we have the following results:
$$\begin{array}{ll} \displaystyle \frac{\partial X}{\partial
w_i}\Big|_{z=0}=\overline w_i\ , & \displaystyle
\frac{\partial X}{\partial \overline w_j}\Big |_{z=0}=w_j\ ,\\
\\
\displaystyle \frac{\partial X}{\partial z_{\alpha}}\Big
|_{z=0}=\frac{\partial X}{\partial z_{kl}}\Big |_{z=0}=0\ ,&\displaystyle
\frac{\partial X}{\partial \overline z_{\beta}}\Big |_{z=0}=
\displaystyle\frac{\partial X}{\partial \overline z_{st}}\Big |_{z=0}=0\ ,\\
\\
\displaystyle \frac{\partial^2 X}{\partial z_{\alpha}\partial
\overline z_{\beta}}\bigg|_{z=0}
=\frac{X}{K}\delta_{ks}\delta_{lt}\ ,& \displaystyle \frac{\partial^2 X}{\partial
z_{\alpha}\partial \overline w_j}\bigg|_{z=0}=0\ ,\\
\\
\displaystyle \frac{\partial^2 X}{\partial w_i \partial \overline
z_{\beta}}\bigg|_{z=0}=0\ ,& \displaystyle \frac{\partial^2 X}{\partial w_i \partial \overline
w_j}\bigg|_{z=0}=\delta_{ij}\ .\\
\end{array}
$$
where
$
\begin{array}{lll}
\delta_{mn}= \left\{
\begin{array}{ll}
1, & m=n\\ 0, & m \not= n\\
\end{array}
\right. \\
\end{array}
$. Applying the above result to formula (3), we obtain
$$\begin{array}{ll}
\displaystyle \frac{\partial^2 g}{\partial z_{\alpha}\partial
\overline z_{\beta}}\bigg|_{z=0} =\displaystyle
\frac{Y}{K}\delta_{\alpha\beta}\ ,& \displaystyle \frac{\partial^2
g}{\partial z_{\alpha}\partial \overline
w_j}\bigg|_{z=0}=\displaystyle \frac{\partial^2
g}{\partial w_i \partial \overline z_{\beta}}\bigg|_{z=0}=0\ ,\\
\\
\displaystyle \frac{\partial^2 g}{\partial w_i \partial \overline
w_j}\bigg|_{z=0}=\displaystyle Y\delta_{ij}+ \displaystyle Y^2
\overline w_{i} w_{j}\ .&
\end{array}$$
Therefore
$$\left(g_{\alpha \overline\beta}(0,w^*)\right)=\left(\begin{array}{cc}
   \frac{Y}{K}I^{(N-r)}&0\\
     0& YI^{(r)} + Y^{2}{\overline w^{*}}^{t}w^{*}
    \end {array}
\right).$$ where $I^{(r)}$ denotes $r\times r$ unit matrix and we
write $w^*$ back instead of $w$.

Notice that $w^{*}= w\det(I-Z\overline{Z})^{-\frac{1}{2K}}$ when
$Z_{0}=Z$, we get
$$\begin{array}{lll}\det (g_{\alpha\beta}(0,w^{*}))&=&
 (\frac{Y}{K})^{N-r}\det(YI^{(r)} + Y^{2}{\overline w^{*}}^{t}w^{*})\\
&=&(\frac{Y}{K})^{N-r}Y^{r}\det\displaystyle(I^{(r)}+{XY}{|w|^{-2}}\overline{w}^{t}w)\\
&=&(\frac{Y}{K})^{N-r}Y^{r}\det\displaystyle(I^{(1)}+{XY}{|w|^{-2}}w\overline {w}^{t})\\
&=& K^{r-N}Y^{N}(1+XY)\\
&=&K^{r-N}(1-X)^{-(N+1)}.\\
     \end {array}$$
Hence
$$\begin{array}{lll} \det(g_{\alpha \overline
\beta}(z,w))&=&  |\mathrm{J}_{F_0}|^2\det(g_{\alpha \overline\beta}(0,w^*))\\
&=&K^{r-N}(1-X)^{-(N+1)}\det(I-Z\overline{Z})^{-(p+1+\frac{r}{K})},
\end{array}\eqno{(4)}$$
while the other side of Monge-Amp\` ere equation is
$$\begin{array}{lll} e^{(N+1) g(z,w)}&=& e^{(N+1)\log [\frac{1}{1-X}\det(I-Z\overline
Z)^{-\frac{1}{K}}K^{\frac{r-N}{1+N}}]}\\
&=&K^{r-N}(1-X)^{-(N+1)}\det(I-Z\overline{Z})^{-\frac{N+1}{K}}.
\end{array}\eqno{(5)}$$
It is easy to obtain that formula (4) and (5) is equal in the case
$K=\frac p2+\frac 1{p+1}$. That is the function $g$ we guess is a
solution to complex Monge-Amp\` ere equation. It remains to prove
that $g$ satisfies the Dirichlet boundary condition.

If $(\stackrel{\sim}{z}, \stackrel{\sim}{w})\in
\partial Y_{II}$ and $\stackrel{\sim}{w}\not=0$, when $(z,w)\in Y_{II}$ and $
(z,w)\rightarrow (\stackrel{\sim}{z}, \stackrel{\sim}{w})$,
 we have $X\rightarrow 1^-$, so $ \frac{1}{1-X} \rightarrow
+\infty$, meanwhile ${\rm{det}}(I-Z\overline Z) \rightarrow
|\stackrel{\sim}{w}|^{2K}>0$. Hence we have $g(z,w) \rightarrow
+\infty,$ as $(z,w)\rightarrow \partial Y_{II}. $

If $(\stackrel{\sim}{z}, \stackrel{\sim}{w})\in \partial Y_{II}$
and $\stackrel{\sim}{w}=0$, when $(z,w)\in Y_{II}$ and
$(z,w)\rightarrow(\stackrel{\sim}{z},0)$, we have$
\frac{1}{1-X}>1$, ${\rm{det}}(I-Z\overline Z) \rightarrow 0,
{\rm{det}}(I-Z\overline Z)^{-\frac{1}{K}} \rightarrow +\infty$, we
also have $g(z,w) \rightarrow +\infty,$ as$(z,w) \rightarrow
\partial Y_{II}.$

Up to now, we have proved our conjecture, that is the function
$$g=\log
[\frac{1}{1-X}\det(I-Z\overline
Z)^{-\frac{1}{K}}K^{\frac{r-N}{1+N}}]$$ generates a complete
Einstein-K$\ddot{\rm a}$hler metric on $Y_{II}(r,p;K)$ in the case
$K=\frac p2+\frac 1{p+1}$. In general, $Y_{II}(r,p;K)$ is a
nonhomogeneous domain when $ K=\frac p2 +\frac 1{p+1}$ and $p>1$.

\section{Holomorphic sectional curvature}

Since in the case $K=\frac p2+\frac 1{p+1}$ the complete
Einstein-K\"{a}hler metric on $Y_{II}(r,p;K)$ is generated by
$$
g=\log[\frac{1}{1-X} \det(I-Z\overline Z)^{-\frac
1K}K^{\frac{r-N}{1+N}}],\qquad N={p(p+1)}/{2}+r,
$$
the holomorphic sectional curvature $\omega[(z,w),d(z,w)]$ on
$Y_{II}(r,p;K)$ under this metric has the following form:
$$
\omega[(z,w),d(z,w)]=\frac{d(z,w)[-\overline{d}dT+dT
T^{-1}\overline{dT}^t]\overline{d(z,w)}^t}
{[d(z,w)T\overline{d(z,w)}^t]^2}. $$ where $$
 d=\sum \frac{\partial }{\partial z_{\alpha}}dz_{\alpha},\quad
\overline d=\sum \frac{\partial }{\partial \overline z_{\alpha}}
\overline {dz}_{\alpha},\quad \alpha=1,2,\cdots,N,\quad
T=\left(\frac{\partial^{2} g} {\partial
z_{\alpha}\partial\overline z_{\beta}}\right) _{1 \leq
\alpha,\beta \leq N}.$$

Now that holomorphic sectional curvature $\omega[(z,w),d(z,w)]$ is
invariant under the holomorphic automorphism group
$\mathrm{Aut}(Y_{II})$, and due to the Lemma 1, for $\forall
(z,w)\in Y_{II}$ there exists $F\in {\rm{Aut}}(Y_{II})$ such that
$F(z,w)=(0,w^*)$. So it suffices to calculate the
$\omega[(z,w),d(z,w)]$ on point $(0,w^*)$. By sec.2, the complete
Einstein-K\"{a}hler metric matrix is $T=\left(
\begin{array}{ll}
T_{11} & T_{12}\\
T_{21} & T_{22}
\end{array}\right)$, where $$\begin{array}{l}
 T_{11}={K^{-1}Y}[A^t\overline A \cdot\!\!\times A^t\overline A]_s
+{X}{K^{-2}Y^2}E(Z)^t\overline{E(Z)},\\
 T_{12}={K^{-1}Y^2}\det(I-Z\overline Z)^{-\frac{1}{K}}E(Z)^tw,\\
 T_{21}=\overline T^t_{12},\\
 T_{22}={Y^2}{\overline{w}^t w} \det(I-Z\overline
Z)^{-\frac{2}{K}}+YI^{(r)}\det(I-Z\overline Z)^{-\frac{1}{K}}.
\end{array}
$$
$Y=(1-X)^{-1}$ and $A,[A^t\overline A \cdot\!\!\times A^t\overline
A]_s, E(Z)$ are the same as before. Note that
$$dT=\left(
\begin{array}{cc}
dT_{11}&dT_{12}\\dT_{21}&dT_{22}\\
\end{array}
\right), \quad \overline{d}dT=\left(
\begin{array}{cc}
\overline{d}dT_{11}&\overline{d}dT_{12}\\
\overline{d}dT_{21}&\overline{d}dT_{22}\\
\end{array}
\right) .$$ Using some known results$^{[12]}$
$$E(Z)|_{z=0}=dE(Z)|_{z=0}=0,\ \ \overline
dE(Z)^t|_{z=0}=(\overline{dz})^t,\ \ d\overline {E(Z)}|_{z=0}=dz$$
and $d[A^t\overline A \cdot\!\!\times A^t\overline A]_s|_{z=0}=0$,
we can obtain $$
\begin{array}{ll}
 dT_{11}|_{z=0}={K}^{-1}Y^2\overline w dw^tI^{(N-r)},
& dT_{12}|_{z=0}=0,\\
 dT_{21}|_{z=0}={K}^{-1}Y^2\overline w^t dz,
& dT_{22}|_{z=0}=\left(2Y^3\overline
w^tw+Y^2I^{(r)}\right)\overline wdw^t+Y^2\overline w^tdw.
\end{array}
$$
furthermore
$$
\begin{array}{lll}
\overline{d}dT_{11}|_{z=0}&=&K^{-1}Y\bigl(Y|dw|^{2}I+2Y^2|w\overline{dw}^{t}|^{2}I+
{K^{-1}}XY|dz|^2I\\

&&+{K^{-1}}XY\overline{dz}^tdz
+\overline d d[A^t\overline A \cdot\!\!\times A^t\overline A]_s|_{Z=0}\bigr),\\

\overline{d}dT_{12}|_{z=0}&=&K^{-1}Y^2\bigl(2Y\overline{w}dw^{t}\overline{dz}^tw+
\overline{dz}^tdw\bigr),\\

\displaystyle \overline{d}dT_{21}|_{z=0}&=&K^{-1}Y^2\bigl(2Y
w\overline{dw}^{t}\overline{w}^{t}dz+ \overline{dw}^{t}dz\bigr),\\

\overline{d}dT_{22}|_{z=0}&=&Y^2\bigl({K}^{-1}(2Y\overline
w^tw+I)|dz|^2
+2Y(3Y\overline w^tw+I)|w\overline {dw}^t|^2+\overline{dw}^tdw\\
&&+(2Y\overline w^tw+I)|dw|^2+2Y(\overline wdw^t\overline{dw}^tw
+w\overline{dw}^t\overline w^tdw)\bigr).\\
\end{array}
$$
If one denote
$$ [-\overline{d}dT+dTT^{-1}\overline{dT}^t]_{z=0} =\left(
\begin{array}{cc}
R_{11}&R_{12}\\ R_{21}&R_{22}
\end{array}
\right),
$$
notice that $$\displaystyle T^{-1}|_{z=0}= \left(
\begin{array}{cc}
\frac{K}{Y}{I}&0\\ 0&\displaystyle
{Y}^{-1}\left(I-\overline{w}^t w\right)\\
\end{array}
\right),
$$
then we get
$$
\begin{array}{lll}
R_{11}&=&-K^{-1}Y\bigl(Y^2|w\overline{dw}^t|^2{I}+K^{-1}XY(\overline{dz}^t
dz+|dz|^2{I})+Y|dw|^2{I}\\

&&+\overline d d[A^t\overline A\cdot\!\!\times A^t\overline A]_{s}|_{Z=0}\bigr),\\

R_{12}&=&-K^{-1}Y^2\bigl(Y\overline w dw^t
\overline{dz}^t w+\overline{dz}^t dw\bigr),\\

R_{21}&=&\overline R_{12}^t,\\

R_{22}&=&-Y^2\bigl({K}^{-1}Y|dz|^2\overline w^t w+
{K}^{-1}|dz|^2{I}+|dw|^2{I}+\overline{dw}^tdw\\

&&+{Y}(|\overline wdw^t|^2{I}+
w\overline{dw}^t \overline{w}^t dw+\overline
wdw^t\overline{dw}^t w+|dw|^2\overline w^t w)\\

&&+2Y^{2}|\overline wdw^t|^2\overline w^t w\bigr).
\end{array}\eqno{(6)}
$$
Complying with formula (6) and the result  $$ dz\overline
dd[A^t\overline A \cdot\!\!\times A^t\overline
A]_s|_{Z=0}\overline {dz}^t = 2tr(dZ\overline
{dZ}dZ\overline{dZ})$$ in ref.[12], one can easily get
$$\begin{array}{lll}
&&(dz,dw)[-\overline{d}dT+dTT^{-1}\overline{dT}^t]\overline{(dz,dw)}^t|_{z=0}\\

&=&dzR_{11}\overline{dz}^t+dwR_{21}\overline{dz}^t+dzR_{12}\overline{dw}^t
+dwR_{22}\overline{dw}^t\\

&=&-{2}{K^{-2}}XY^2|dz|^4-4{K}^{-1}Y^3|dz|^2|\overline{w}dw^t|^2-2Y^2|dw|^4-2Y^4|\overline wdw^t|^4\\

&&-{4}{K}^{-1}Y^2|dw|^2|dz|^2-4Y^3|\overline{w}dw^t|^2|dw|^2
-{2K^{-1}Y}\textrm{tr}(dZ\overline {dZ} dZ\overline {dZ}).
\end{array}
$$
Additionally,
$$(dz,dw)T\overline{(dz,dw)}^t|_{z=0}={K}^{-1}{Y}|dz|^2+Y^2|\overline{w}dw^t|^2+Y|dw|^2,$$
thus the holomorphic sectional curvature on point$(w,0)$ under the
complete Einstein-K\"{a}hler metric is
$$\omega((z,w),d(z,w))|_{z=0}
=-2+\frac{ 2{K^{-2}}Y|dz|^4-{2{K}^{-1}Y}{tr}(dZ\overline
{dZ}dZ\overline {dZ})}
{({K}^{-1}{Y}|dz|^2+Y^2|\overline{w}dw^t|^2+Y|dw|^2)^2}.$$

It is apparent that if let $p=1$, then $K=1, {tr}(dZ\overline
{dZ}dZ\overline {dZ})=|dz|^4.$ In this case, $Y_{II}(r,p;K)$ is a
unit ball in $\mathbb{C}^{r+1}$ and
$\omega((z,w),d(z,w))|_{z=0}=-2$. It is a well-known result. Our
work is to give its estimate in the case $p>1$.

Since $|dz|^2=tr(dZ\overline{dZ})$, according to lemma 5, we get
$${p}^{-1}|dz|^4 \leq tr(dZ\overline {dZ} dZ\overline {dZ}) \leq
|dz|^4.$$ Applying it to $\omega((z,w),d(z,w))|_{z=0}$, we have
$$2(1-K)\frac{Y}{K^2}|dz|^4\leq
\frac{2Y}{K^{2}}|dz|^4-{\frac{2Y}{K}}{tr}(dZ\overline
{dZ}dZ\overline {dZ}) \leq 2(1-\frac{K}p)\frac{Y}{K^2}|dz|^4.
$$
Now that $1-K<0$ and $1-\frac{K}p>0$ hold in the case $p>1$ and
notice that $Y\geq 1 (0\leq X < 1)$, the above inequality can be
expanded as
$$\begin{array}{lll}
&&2(1-K){({K}^{-1}{Y}|dz|^2+Y^2|\overline{w}dw^t|^2+Y|dw|^2)^2}\\
&\leq&2{K^{-2}}Y|dz|^4-{2{K}^{-1}Y}{tr}(dZ\overline {dZ}dZ\overline {dZ})\\
&\leq&
2(1-{K}p^{-1}){({K}^{-1}{Y}|dz|^2+Y^2|\overline{w}dw^t|^2+Y|dw|^2)^2}.
\end{array}
$$
Applying it in the expression of $\omega((z,w),d(z,w))|_{z=0}$ one
can obtain the following result immediately:
$$
-2K\leq \omega((z,w),d(z,w))\leq -\frac{2K}p
$$
where $K=\frac p2+\frac 1{p+1}$ and $p>1$. This estimate is the
sharp estimate because of the following facts:

At the point $(z,0)$ and direction $(dz,0)$  one easily knows that
$X=0$ and $Y=1$, and
$$\begin{array}{lll}
\omega((z,w),d(z,w))|_{z=0} =\displaystyle
-2+\frac{2{K^{-2}}|dz|^4-{2{K}^{-1}}{tr}(dZ\overline {dZ}dZ\overline
{dZ})}{{K}^{-2}|dz|^4}\ .
\end{array}
$$
If choose
$$
dZ=U^t\left(
\begin{array}{ccc}
\lambda_1&&{0}
\\ {0}&&{0}\\
\end{array}
\right ) U ,\quad  \lambda_1>0,$$ then according to lemma 5 one has
$tr(Z\overline{Z}Z\overline{Z})= tr(Z\overline{Z})
tr(Z\overline{Z})=|dz|^4$, which implies
$\omega((z,w),d(z,w))|_{z=0}=-2K.$ \\
If choose
$$ dZ=U^t\left(
\begin{array}{ccc}
\lambda_1& &{0} \\ &\ddots&\\ {0}&&\lambda_p\\
\end{array}
\right )  U ,\quad \lambda_1= \lambda_2= \cdots= \lambda_p>0,
$$
then one has $tr(Z\overline{Z}Z\overline{Z})=p^{-1}
tr(Z\overline{Z}) tr(Z\overline{Z})=p^{-1}|dz|^4$ by lemma 5, which
implies $\omega((z,w),d(z,w))|_{z=0}=-\frac{2K}p.$

Therefore our estimate is the sharp estimate.

\section{K\"{a}hler--Einstein metric Is Equivalent to the Bergman metric}

In this section we will prove that the K\"{a}hler--Einstein metric
is equivalent to the Bergman metric on $Y_{II}(r,p;K)$ in the case
$K=\frac p2+\frac 1{p+1}$.
\\
\\
{\bf Definition:} Let $\mathcal{B}_0$ and $\mathcal{E}_0$ be two
complete metrics on domain $\Omega$. If there exists two positive
constant $a$ and $b$(with $a\geq b$) such that the following
inequality holds:
$$
b\leq \frac{\mathcal{B}_0}{\mathcal{E}_0}\leq a.
$$
Then we called that the metric $\mathcal{B}_0$  is equivalent to the
metric $\mathcal{E}_0.$
\\

From sec.2 we know in the case $K=\frac p2+\frac 1{p+1}$ the
complete K\"{a}hler--Einstein metric on $Y_{II}(r,p;K)$ is
$$
\mathcal{E}_{II}:=(dz,dw)T_{\mathcal{E}_{II}}\overline{(dz,dw)}^t,
$$
where
$$
T_{\mathcal{E}_{II}}=\mathrm{J}_{F_0}T^{(0)}_{\mathcal{E}_{II}}
\overline {\mathrm{J}}_{F_0}^t=\mathrm{J}_{F_0}
\left(\begin{array}{cc}
   \frac{Y}{K}I^{(N-r)}&0\\
     0& YI^{(r)} + Y^{2}{\overline w^{*}}^{t}w^{*}
    \end {array}
\right) \overline {\mathrm{J}}_{F_0}^t.$$

From ref.[11], we know the Bergman kernel function on
$Y_{II}(r,p;K)$ is
$$
K_{II}(W,Z;\overline{W},\overline{Z})=K^{-\frac{p(p+1)}{2}}\pi^{-\frac{p(p+1)}{2}-r}
G(Y)\mathrm{det}(I-Z\overline{Z})^{-(p+1+\frac rK)}
$$
where
$$
G(Y)=\sum^h_{j=0}b_j\Gamma(r+j)Y^{r+j},\quad
h=\frac{p(p+1)}{2}+1,\quad b_h=2^{\frac{p(p+1)}{2}+1},\quad
\Gamma(r+j)=(r+j-1)!.
$$
Thus the Bergman metric on $Y_{II}(r,p;K)$ has the following form:
$$
\mathcal{B}_{II}:=(dz,dw)T_{\mathcal{B}_{II}}\overline{(dz,dw)}^t,
$$
where the metric matrix
$$
T_{\mathcal{B}_{II}}=\left(\frac{\partial^{2} \mathrm{log}K_{II}}
{\partial z_{\alpha}\partial\overline z_{\beta}}\right)=
\mathrm{J}_{F_0}\left(\frac{\partial^{2} \mathrm{log}K_{II}}
{\partial z_{\alpha}\partial\overline
z_{\beta}}\right)_{Z_0=Z}\overline {\mathrm{J}}_{F_0}^t
=\mathrm{J}_{F_0}T^{(0)}_{\mathcal{B}_{II}} \overline
{\mathrm{J}}_{F_0}^t$$ and
$$ T^{(0)}_{\mathcal{B}_{II}}=
\left(\begin{array}{cc}
   (\frac{1}{K}H'X+p+1+\frac rK)I^{(N-r)}&0\\
     0& H'I^{(r)} + H''{\overline w^{*}}^{t}w^{*}
    \end {array}
\right),\quad H=\mathrm{log}G(Y).
$$

If we denote $(dz,dw)\mathrm{J}_{F_0}$ by
$(d\mathcal{Z},d\mathcal{W})$, then $\mathcal{E}_{II}$ and
$\mathcal{B}_{II}$ can be rewritten as
$$\begin{array}{lll}
\mathcal{E}_{II}&=&\frac YK \left|d\mathcal{Z} \right|^2+
d\mathcal{W}(YI^{(r)} + Y^{2}{\overline
w^{*}}^{t}w^{*})\overline{d\mathcal{W}}^t
\\
\mathcal{B}_{II}&=&(\frac{1}{K}H'X+p+1+\frac rK)\left|d\mathcal{Z}
\right|^2+ d\mathcal{W}(H'I^{(r)} + H''{\overline
w^{*}}^{t}w^{*})\overline{d\mathcal{W}}^t.
\end{array}
$$
According to ref.[13] one knows that the vector
$w^*=(w^*_1,w^*_2,\cdots,w^*_r)$ can be transformed into
$$w^*=e^{i\theta}(\lambda,0,\cdots,0)U,\quad \lambda\geq¡¡0,$$
where $U$ is a unitary matrix. Hence $$
\overline{w^*}^tw^*=\overline{U}^t\left(\begin{array}{cc}\lambda^2&0\\
0&\mathbf{0} \end{array}\right)U,\quad H'I^{(r)} + H''{\overline
w^{*}}^{t}w^{*}=\overline{U}^t\left(\begin{array}{cc}H'+H''\lambda^2&0\\
0&H'I^{(r-1)} \end{array}\right)U$$ and
$$
YI^{(r)} + Y^{2}{\overline w^{*}}^{t}w^{*}=
\overline{U}^t\left(\begin{array}{cc}Y+Y^2\lambda^2&0\\
0&YI^{(r-1)} \end{array}\right)U.
$$
Let $d\mathcal{W}\overline{U}^t=(d\mathbf{w},d\mathbf{W})$, then
$\mathcal{E}_{II}$ and $\mathcal{B}_{II}$ have the following form:
$$\begin{array}{lll}
\mathcal{E}_{II}&=&\frac{Y}{K}\left|d\mathcal{Z}
\right|^2+(Y+Y^2\lambda^2)\left|d\mathbf{w}\right|^2+Y\left|d\mathbf{W}\right|^2,\\
\mathcal{B}_{II}&=&(\frac{1}{K}H'X+p+1+\frac rK)\left|d\mathcal{Z}
\right|^2+(H'+H''\lambda^2)\left|d\mathbf{w}\right|^2+H'\left|d\mathbf{W}\right|^2.
\end{array}$$
It is known $T_{\mathcal{E}_{II}}$ and $T_{\mathcal{B}_{II}}$ are
positive definite matrixes, which implies
$$\frac{1}{K}H'X+p+1+\frac rK>0,\quad  H'+H''\lambda^2>0, \quad  H'>0,\quad Y>0.$$
Therefore, if we denote
$$
\Phi(X)=\frac{{K}^{-1}H'X+p+1+rK^{-1}}{{Y}{K}^{-1}},\quad
\Psi(X)=\frac{H'+H''\lambda^2}{Y+Y^2\lambda^2},\quad
\Upsilon(X)=\frac{H'}{Y},
$$
then $\Phi(X),\Psi(X)$ and $\Upsilon(X)$ are all positive continues
functions of $X$ on the interval $[0,1)$. If $$ \lim_{X\rightarrow
1}\Phi(X), \lim_{X\rightarrow 1}\Psi(X),
\lim_{X\rightarrow1}\Upsilon(X)$$ are existent and positive, then
all of $\Phi(X), \Psi(X), \Upsilon(X)$ have the positive maximum and
the positive minimum on [0, 1) respectively.

We know that
$$
G(Y)=\sum^h_{j=0}b_j\Gamma(r+j)Y^{r+j},\quad
h=\frac{p(p+1)}{2}+1,\quad b_h=2^{\frac{p(p+1)}{2}+1},
$$
then
$$
\frac{dG(Y)}{dX}=G'(Y)=\sum^h_{j=0}b_j\Gamma(r+j+1)Y^{r+j+1},
\frac{d^2G(Y)}{dX^2}=G''(Y)=\sum^h_{j=0}b_j\Gamma(r+j+2)Y^{r+j+2},\eqno{(7)}
$$
and
$$H'=G'(Y)G^{-1}(Y),\quad H''=G''(Y)G^{-1}(Y)-G'(Y)^2G^{-2}(Y).$$

We shall compute the limits of $\Phi(X),\Psi(X)$ and $\Upsilon(X)$
when $X$ tends to 1.
$$
\lim_{X\rightarrow 1}\Phi(X)=\lim_{Y\rightarrow\infty}\Phi(X)=
\lim_{Y\rightarrow\infty}\frac{{K}^{-1}H'X+p+1+rK^{-1}}{{Y}{K}^{-1}}
=\lim_{Y\rightarrow\infty}\frac{G'(Y)}{G(Y)Y},
$$
applying formula (7), we have
$$
\lim_{X\rightarrow 1}\Phi(X)=\lim_{Y\rightarrow\infty}\Phi(X)=
\frac{b_h\Gamma(r+h+1)Y^{r+h+1}}{b_h\Gamma(r+h)Y^{r+h+1}}=r+h=N+1,
$$
where $N=\frac{p(p+1)}2+r$ is the dimension of $Y_{II}(r,p;K)$.

Therefore there exists $0<\nu<\mu$ such that $$0<\nu\leq \Phi(X)\leq
\mu. $$

Similarly, we can also obtain
$$
\lim_{X\rightarrow 1}\Psi(X)=\lim_{Y\rightarrow\infty}\Psi(X)=
\lim_{Y\rightarrow\infty}\frac{H'+H''\lambda^2}{Y+Y^2\lambda^2}
=\lim_{Y\rightarrow\infty}\frac{H''}{Y^2}=N+1
$$
and
$$
\lim_{X\rightarrow1}\Upsilon(X)=\lim_{Y\rightarrow\infty}\Upsilon(X)
=\lim_{Y\rightarrow\infty}\frac{H'}{Y}=N+1.
$$
Therefore there also exist positive constants $\zeta, \eta, \rho$
and $\varrho$ such that
$$
\zeta\leq\Psi(X)\leq\eta,\quad\rho\leq\Upsilon(X)\leq\varrho. $$

Let $a=max\{\mu,\eta,\varrho\}$ and $b=min\{\nu,\zeta,\rho\}$,
then we have
$$
b\leq \frac{\mathcal{B}_{II}}{\mathcal{E}_{II}}\leq a.
$$
Up to now we can say that the complete K\"{a}hler--Einstein metric
is equivalent to the Bergman metric on $Y_{II}(r,p;K)$ in the case
$K=\frac p2+\frac 1{p+1}$.

\end{document}